\documentclass{article}
\usepackage{graphicx} % Required for inserting images
\usepackage{amsfonts}
\usepackage{amscd,color}
\usepackage{amsmath,amsfonts,amssymb,amscd}
\usepackage{booktabs}
\usepackage{mathrsfs}
\usepackage{indentfirst,graphicx,epsfig}
\usepackage{xcolor}
\usepackage{epstopdf}
\usepackage{caption}
\usepackage{ifpdf}
\usepackage{float}
\usepackage{footnote}
\usepackage{epsfig}
\usepackage{latexsym}
\usepackage{caption}
\usepackage{amsmath,amsthm}
\usepackage{amsfonts}
\usepackage{amssymb}
\usepackage{graphicx}
\usepackage{graphics}
\usepackage{bm}
\usepackage{color}
\usepackage{multirow}
\usepackage{subfigure}
\usepackage{threeparttable}
\usepackage{ulem}

\makeatletter

\newcommand{\Rmnum}[1]{\expandafter\@slowromancap\romannumeral #1@}
\makeatother
\makeatletter
\let\@fnsymbol\@arabic
\makeatother

\makeatletter
  \newcommand\figcaption{\def\@captype{figure}\caption}
  \newcommand\tabcaption{\def\@captype{table}\caption}
\makeatother
\begin{document}
\newtheorem{theorem}{Theorem}[section]
\newtheorem{observation}[theorem]{Observation}
\newtheorem{corollary}[theorem]{Corollary}
\newtheorem{algorithm}[theorem]{Algorithm}
\newtheorem{problem}[theorem]{Problem}
\newtheorem{question}[theorem]{Question}
\newtheorem{lemma}[theorem]{Lemma}
\newtheorem{proposition}[theorem]{Proposition}
\newtheorem{definition}[theorem]{Definition}
\newtheorem{guess}[theorem]{Conjecture}
\newtheorem{claim}[theorem]{Claim}
\newtheorem{example}[theorem]{Example}
\newtheorem{acknowledgement}[theorem]{Acknowledgement}
\newtheorem{axiom}[theorem]{Axiom}
\newtheorem{case}[theorem]{Case}
\newtheorem{conclusion}[theorem]{Conclusion}
\newtheorem{conjecture}[theorem]{Conjecture}
\newtheorem{criterion}[theorem]{Criterion}
\newtheorem{exercise}[theorem]{Exercise}
\newtheorem{notation}[theorem]{Notation}
\newtheorem{solution}[theorem]{Solution}
\newtheorem{summary}[theorem]{Summary}
\newtheorem{fact}[theorem]{Fact}
\newtheorem{remark}[theorem]{Remark}
\newtheorem{subclaim}[theorem]{Subclaim}

\newcommand{\pp}{{\it p.}}
\newcommand{\de}{\em}
\newcommand{\mad}{\rm mad}

\newcommand*{\QEDA}{\hfill\ensuremath{\blacksquare}}
\newcommand*{\QEDB}{\hfill\ensuremath{\square}}
\newcommand{\xz}[1]{{\color{purple}[XZ: #1]}}
\newcommand{\zw}[1]{{\color{blue} [zw: #1]}}
\newcommand{\qf}{Q({\cal F},s)}
\newcommand{\qff}{Q({\cal F}',s)}
\newcommand{\qfff}{Q({\cal F}'',s)}
\newcommand{\f}{{\cal F}}
\newcommand{\ff}{{\cal F}'}
\newcommand{\fff}{{\cal F}''}
\newcommand{\fs}{{\cal F},s}
\newcommand{\g}{\gamma}
\newcommand{\wrt}{with respect to }

\title{Finiteness of Non-decomposable Critically Frustrated Signed Graphs}
\author{\normalsize Zhiqian Wang$^{1}$\\
\small $^{1}$School of Mathematical Sciences, Zhejiang Normal University, Jinhua, China\\
\small Email: 1522686578@qq.com\\
}

\date{}
\maketitle

\begin{abstract}
A signed graph $(G,\sigma)$ is a graph $G$ with a signature $\sigma$ labeling each edge with a positive or negative sign. Two signatures of $G$ are switching equivalent if one is obtained from the other by sequentially changing the signs of all edges in an edge-cut. The frustration index of a signed graph $(G, \sigma)$ is the minimum number of negative edges among all signatures switching equivalent to $\sigma$. A signed graph is critically $k$-frustrated if it has frustration index $k$, and the removal of any edge decreases its frustration index. A critically $k$-frustrated signed graph is prime if it is irreducible and contains no two edge-disjoint
negative cycles. Steffen and Naserasr et al. asked if there are only finitely many prime critically $k$-frustrated signed graphs for any positive integer $k$. The cases $k\leq 3$ were proved to be true by Cappello et al., and later, the cases $k=4$ and $5$ were verified by Z. Wang. In this paper, the author proves the conjecture by a topological argument.
\end{abstract}

\section[Introduction]{Introduction}
A \textit{signed graph} $(G,\sigma)$ is a graph $G$ together with a \textit{signature} $\sigma$, which is a map $\sigma$: $E(G)\rightarrow \{+,-\}$ that assigns to each edge a positive sign or a negative sign. An edge $e$ of $(G,\sigma)$ is called \textit{negative} if $\sigma(e)=-$ and \textit{positive} otherwise. The sets of positive and negative edges are denoted by $E_\sigma^+(G)$ and $E_\sigma^-(G)$, respectively. We denote by $(G,-)$ (or $(G,+)$) the signed graph with all edges negative (or positive).

For $X\subseteq V(G)$, let $X^c:=V(G)\setminus X$. Let $\partial_G(X):=\{xy\in E(G): x\in X, y\in X^c\}$, which is called an \textit{edge cut} of $G$. The cardinality of $\partial_G(X)$ is denoted by $d_G(X)$. 
 For a signed graph $(G, \sigma)$, let 
\begin{displaymath}
    d_{(G,\sigma)}^+(X)=|\partial_G(X)\cap E_\sigma^+(G)| \mbox{ and }  d_{(G,\sigma)}^-(X)=|\partial_G(X)\cap E_\sigma^-(G)|.
\end{displaymath}

An edge cut $\partial_G(X)$ is called an $(a,b)$-cut if $d_{(G,\sigma)}^+(X)=a$ and $d_{(G,\sigma)}^-(X)=b$. An $(a,a)$-cut is called an \textit{equilibrated cut}.

A \textit{cycle} in $G$ is a connected 2-regular subgraph. A cycle is   \textit{positive} (or \textit{negative}) if it contains an even number (or odd number) of negative edges. A signed graph $(G, \sigma)$ is called \textit{balanced} if it contains no negative cycle and \textit{unbalanced} otherwise.

For $X\subseteq V(G)$, \textit{switching} at $X$ means reversing the signs of the edges in $\partial_G(X)$ and leaving all other signs unchanged. Two signatures $\sigma$ and $\pi$ are \textit{switching equivalent} if one is obtained from the other by a series of switching at some vertex set. We say that $(G,\sigma)$ and $(G,\pi)$ are switching equivalent when $\sigma$ and $\pi$ are switching equivalent signatures. It was shown in \cite{82} that two signatures on the same graph are switching equivalent if and only if they have the same set of negative cycles.

The \textit{frustration index} of a signed graph $(G,\sigma)$ is a fundamental parameter measuring the extent to which a signed graph fails to be balanced. Denoted by $l(G,\sigma)$, it is defined as
\begin{displaymath}
l(G,\sigma)=\min\{|E|\subseteq E(G),\ G\setminus E\text{ is balanced}\}.
\end{displaymath}

This concept was formally developed by Harary \cite{harary} in 1959, motivated by the problem of restoring structural balance in networks of positive and negative relationships. Harary \cite{harary} also established that the minimum number of edges whose deletion makes a signed graph balanced equals the minimum number of edges whose signs must be reversed to achieve balance; moreover, the corresponding minimum edge sets coincide. In 1982, the equivalent modern formulation using switching was developed by T. Zaslavsky \cite{82}, inspired by Seidel's work on Seidel matrix of a graph \cite{siedel}. In the language of switching, the definition of frustration index can be expressed as
\begin{displaymath}
l(G,\sigma)=\min\{|E_\pi^-(G)|:(G,\pi) \mbox{ is switching equivalent to } (G,\sigma)\}.
\end{displaymath}
If $l(G,\sigma)=k$, then $(G,\sigma)$ is said to be \textit{$k$-frustrated}. A signature $\sigma$ is said to be a \textit{minimum signature} if $|E_\sigma^-(G)|=l(G,\sigma)$. Thus if $\sigma$ is a minimum signature and $\partial(X)$ is an $(a,b)$-cut of $(G, \sigma)$, then $a \ge b$.

Naturally, one may wonder what the simplest signed graphs with a prescribed frustration index look like. In \cite{22}, the concept of a critically $k$-frustrated signed graph was introduced by E. Steffen et al.
\begin{definition}
A signed graph $(G,\sigma)$ is \textit{critically $k$-frustrated} if $l(G,\sigma)=k$ and for every $e\in E(G)$, $l(G\setminus e,\sigma)=k-1$.
\end{definition}

We call $(G,\sigma)$ critical when its frustration index need not be specified. In \cite{22}, some basic properties of critical signed graphs were proved.
\begin{theorem}\label{cut}
(\cite{22}) Let $k$ be a positive integer and $(G,\sigma)$ be a $k$-frustrated signed graph. The following statements are equivalent:
\begin{enumerate}
    \item $(G,\sigma)$ is critically $k$-frustrated.
    \item For every edge $e\in E(G)$, there exists a signature $\sigma^\prime$ equivalent to $\sigma$  such that $|E_{\sigma^\prime}^-(G)|=k$ and $e\in E_{\sigma^\prime}^-(G)$.
    \item If $|E_{\sigma}^-(G)|=k$, then every $e\in E_\sigma^+(G)$ is contained in an equilibrated edge cut of $(G, \sigma)$.
\end{enumerate}
\end{theorem}

Assume $t \ge 1$ and let $E_{xy}$ be a set of $t$ parallel edges joining $x$ and $y$, all with the same sign. \textit{Subdividing} $E_{xy}$ means deleting $E_{xy}$ and adding a new vertex $w$, as well as a set $E^+_{xw}$ of $t$ positive  parallel edges between $x$ and $w$  and a set $E_{wy}$ of $t$ parallel edges between $w$ and $y$ that have the same sign as $E_{xy}$. A signed graph $(G', 
\sigma')$ is a \textit{subdivision} of $(G, \sigma)$ if $(G', 
\sigma')$ is obtained from $(G, \sigma)$ by sequentially subdividing some multiedges, and furthermore, if $(G', \sigma')\neq (G, \sigma)$, $(G', \sigma')$ is a \textit{proper subdivision} of $(G, \sigma)$. We say a signed graph is \textit{irreducible} if it is not a proper subdivision of any signed graph.

A critically $k$-frustrated signed graph $(G,\sigma)$ is said to be \textit{$(k_1,\ldots,k_t)$ -decomposable} if $E(G)$ has a partition $E_1 \cup \ldots \cup E_t$ with $t \ge 2$ such that for $i\in\{1,\ldots,t\}$, the signed subgraph $(G[E_i],\sigma)$ is critically $k_i$-frustrated and $k=k_1+\ldots + k_t$. We simply say $(G,\sigma)$ is \textit{decomposable} when the parameters $k_1,\ldots,k_t$ need not be specified. A critical signed graph is \textit{indecomposable} if it is not decomposable.

In \cite{22}, the authors studied a special class of irreducible critical signed graphs arising from the multiple weak 2-linkage problem proposed by Y. Lu et al. \cite{link}. The problem asks the following.
\begin{problem}
(Multiple Weak $2$-Linkage Problem) Let $k\geq 2$ be an integer, and $\{x_1,\ldots,x_k, y_1,\ldots,y_k\}$ be $2k$ fixed vertices in a graph $G$. Does there exist a pair of integers $1\leq i<j\leq k$ such that the graph $G$ contains a pair of edge-disjoint paths $P_i$ and $P_j$ such that $P_r$ joins $x_r$ and $y_r$ for each $r=i,j$?
\end{problem}
It was shown in \cite{link} that if we add a negative edge between every pair $\{x_i, y_i\}$ in $G$, with all original edges assigned positive signs, then Multiple Weak $2$-Linkage Problem is equivalent to checking whether there exist two edge-disjoint negative cycles in the resulting signed graph. This motivates the following definition of prime critical signed graphs.
\begin{definition}\label{def-prime}
A critical signed graph $(G, \sigma)$ is \textit{prime} if it is irreducible and contains no pair of edge-disjoint negative cycles.  
\end{definition}
Since every negative cycle is a critically $1$-frustrated signed graph, in terms of decomposability, this gives the following equivalent description.
\begin{observation}
A critical signed graph $(G,\sigma)$ is prime if and only if $(G,\sigma)$ is irreducible and every critical subgraph is indecomposable.
\end{observation}
We use $\mathcal{S}^*(k)$ to denote the family of prime critical signed graphs with frustration index $k$ in this paper. An example of $\mathcal{S}^*(k)$, called Escher Walls, was constructed for every fixed integer $k$ in \cite{22}. Furthermore, in \cite{25}, it was conjectured that $\mathcal{S}^*(k)$ contains only finitely many signed graphs.
\begin{conjecture}
\label{conj-key}
For any positive integer $k$, $\mathcal{S}^*(k)$ is finite. 
\end{conjecture}
This conjecture has been verified for small $k$. In \cite{22}, it was proved that $\mathcal{S}^*(1)=\{C_{-1}\}$ and $\mathcal{S}^*(2)=\{(K_4,-)\}$; the case $k=3$ was confirmed in \cite{25}; later in \cite{mein}, the author showed that the conjecture holds for $k=4$ or $5$.

Derived from the idea in his previous work \cite{mein}, the author has developed a topological approach to this problem in the current paper. Consequently, the author verifies the conjecture with a structural description in which each member in $\mathcal{S}^*(k)$ has a contraction core in a family with finitely many graphs. As a corollary, the size of signed graphs in $\mathcal{S}^*(k)$ satisfies quadratic lower and upper bounds determined only by $k$. These results are distilled into Theorem \ref{essence1}, Theorem \ref{essence2} and Corollary \ref{quad} in the following, and the terminology used in the statements will be introduced later in Sections 2 and 3.

\begin{theorem}\label{essence1}
Let $\mathcal{T}_k$ be the delta-wye equivalence class of the $k\times k$ M\"{o}bius grid. For $(G,\sigma)\in \mathcal{S}^*(k)$, where $\sigma$ is minimum, there exists a subset $E'\subseteq E^+_\sigma(G)$ such that $E'=\emptyset$ or $E'$ induces a forest in $G$, and $G/E'$ is a $k$-frustrated projective signed graph of $\mathcal{T}_k$.
\end{theorem}

\begin{theorem}\label{essence2}
For $k\geq 4$, the following claims hold:
\begin{itemize}
\item every graph in $\mathcal{T}_k$ is uniquely embedded in $\mathbb{RP}^2$;
\item every tree expansion of graphs in $\mathcal{T}_k$ is uniquely embedded in $\mathbb{RP}^2$;
\item For a signed graph $(G,\sigma)$ with $|E^-_\sigma(G)|=k$, $(G,\sigma)\in \mathcal{S}^*(k)$ if and only if $(G,\sigma)$ is a positive tree expansion of a $k$-frustrated projective signed graph of $\mathcal{T}_k$, and for every $e\in E(G)$, $\rho(\Psi(G\setminus e))<k$, where $\Psi$ is the unique embedding of $G$ in $\mathbb{RP}^2$.
\end{itemize}
\end{theorem}

\begin{corollary}\label{quad}
For any integer $k\geq 3$, $2k^2-k\leq|E(G)|\leq 4k^2-2k-3$ for every $(G,\sigma)\in\mathcal{S}^*(k)$.
\end{corollary}

The remainder of the paper is organized as follows. In Section 2, we introduce some basic concepts involved in the solution of Conjecture \ref{conj-key}. Section 3 proves the main results.

\section{Preliminaries}
In a surface $\Sigma$, a homeomorphic image of the interval $[0,1]$ is called a \textit{simple arc}. A homeomorphic image of $S^1$ is called a \textit{simple closed curve}. A simple closed curve which is homotopic to a point is called a \textit{trivial curve}, otherwise it is called an \textit{essential curve}.

An \textit{embedding} $\Gamma$ of a graph $G$ in a surface $\Sigma$ is an injective continuous map from $G$ in $\Sigma$, in which vertices are associated with points, edges are associated with homeomorphic images of $[0,1]$, and faces are connected components of $\Sigma-\Gamma(G)$. The image $\Gamma(G)$ can be regarded as a well defined graph isomorphic to $G$, and the set of vertices, edges and faces of $\Gamma(G)$ are denoted by $V(\Gamma)$, $E(\Gamma)$ and $F(\Gamma)$ respectively. By convention, we omit $\Gamma$ when the embedding is clear from context. Let $\epsilon(G)$ and $\tilde{\epsilon}(G)$ denote the orientable and nonorientable
genera of the connected graph $G$, respectively.

An embedding $\Gamma$ is a \textit{2-cell embedding} if every face is homeomorphic to an open disk. Two embeddings $\Gamma_1$, $\Gamma_2$ are said to be \textit{equivalent} if there exists a self-homeomorphism $\eta:\Sigma\mapsto\Sigma$ such that $\Gamma_2=\eta\circ\Gamma_1$.

\begin{lemma}\label{genus}
Denote the orientable and nonorientable genus by of a connected graph $G$ by $\epsilon(G)$, and $\widetilde{\epsilon}(G)$ respectively, and let $\widetilde{\Sigma}$ be the nonorientable surface of genus $\widetilde{\epsilon}(G)$. Then every embedding of $G$ in $\widetilde{\Sigma}$ is a 2-cell embedding if and only if $\widetilde{\epsilon}(G)\leq 2\epsilon(G)$.
\end{lemma}

A graph $G$ is said to be \textit{uniquely embedded} in the surface $\Sigma$ if every two embeddings of $G$ in $\Sigma$ are equivalent. If both $\Gamma_1$ and $\Gamma_2$ are 2-cell embeddings, and $G$ is a labeled graph(automorphism not considered), then there is another formulation to describe equivalence.
\begin{lemma}\label{2cell}
(\cite{gos})Let $G$ be a signed graph, and $\Gamma_1,\Gamma_2: G\mapsto \Sigma$ be 2-cell embeddings. Then $\Gamma_1$ and $\Gamma_2$ are equivalent if and only if the sets of facial boundary walks of $\Gamma_1(G)$ and $\Gamma_2(G)$ are identical.
\end{lemma}

The \textit{representativity} (or \textit{face width}) of an embedding is a measure of how densely the graph is embedded in the surface. It is defined as follows.
\begin{definition}
Let $\Gamma: G\mapsto \Sigma$ be an embedding. Then the representativity $\rho(\Gamma)$ of $\Gamma$ is 
\begin{displaymath}
\rho(\Gamma)=\min\{|\gamma\cap \Gamma(G)|\,|\gamma\,\, \mbox{is an essential curve in }\Sigma\}.
\end{displaymath}
\end{definition}
We may use $\rho(G)$ instead sometimes when the embedding $\Gamma$ is clear from the context and needs not to be specified. Any essential curve $\gamma$ which intersects $\Gamma(G)$ in exactly $\rho(\Gamma)$ points is called a \textit{minimal essential curve} respect to $\Gamma$.

\begin{figure}[H]
\centering
\includegraphics[width=8cm]{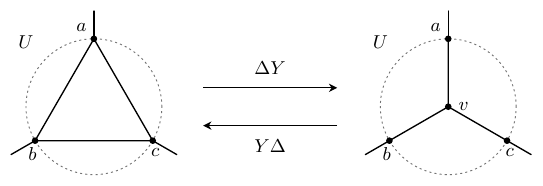}
\caption{$\Delta Y$ and $Y\Delta$ operation.}
\label{delta_y}
\end{figure}

A \textit{triad} of a graph $G$ consists of a degree-three vertex and its three incident edges. Suppose $\Gamma: G\mapsto \Sigma$ is a embedding and $\triangle$ is a triangle in $\Gamma(G)$ bounding a face which is homeomorphic to an open disk. The \textit{$\Delta Y$ operation} at $\triangle$ is to delete the three edges of $\triangle$, and then add a triad centered at a new vertex $v$. Conversely, the \textit{$Y\Delta$ operation} is to replace a triad with a triangle bounding a new face homeomorphic to an open disk. The two reciprocal operations are illustrated in Figure \ref{delta_y}. Both operations are performed in a local disk $U$ containing only the three edges of the triangle and the triad, and the remaining edges not involved in the operations lie in the exterior of $U$.
Given two embeddings $\Gamma_1: G_1\mapsto \Sigma$ and $\Gamma_2: G_2\mapsto \Sigma$, $\Gamma_1(G_1)$ and $\Gamma_2(G_2)$ are \textit{delta-wye equivalent} if there exists a series (may be empty) of $\Delta Y$ and $Y\Delta$ operations that transforms $\Gamma_1(G_1)$ into $\Gamma_2(G_2)$. It is obvious that delta-wye equivalence is an equivalence relation.

\section{Finiteness of $\mathcal{S}^*(k)$}
\subsection{Frustration index and representativity}
To apply the topological notions introduced in Section 2, we start from a description of signed graphs in $\mathcal{S}^*(k)$ \cite{22}, which is derived from the characterization provided by \cite{link}, of graphs for which the multiple weak 2-linkage problem has a negative answer.
\begin{theorem}\label{basis}
Let $k\geq 2$ and $(G,\sigma)$ be a signed graph in $\mathcal{S}^*(k)$ with $E^-_\sigma(G)=\{x_1y_1,\ldots, x_ky_k\}$. Then there exists a permutation $\tau$ on $[k]$, and $(G,\sigma)$ is obtained from a 2-connected plane cubic graph $G'$ (all its edges are positive) by selecting a facial circuit $C$, inserting $2k$ distinct vertices $z_1,\ldots, z_{2k}$ on the edges of $C$ in this cyclic order, and adding $k$ negative edges $e^-_i=x_iy_i$, where $\{z_i, z_{k+i}\}=\{x_{\tau(i)}, y_{\tau(i)}\}$ for $i\in [k]$.

Furthermore, for $k\geq 3$, if $(G,\sigma)\in \mathcal{S}^*(k)$, then the underlying graph $G$ is essentially 4-edge-connected.
\end{theorem}
The projective plane can be drawn in the form of a plane disk with an outer dotted circle to emphasize the border of the crosscap. By Theorem \ref{basis}, a signed graph $(G, \sigma) \in \mathcal{S}^*(k)$ admits a drawing consisting of a plane
graph together with $k$ negative edges drawn through the crosscap. Such a drawing of $(G, \sigma)$ is called a \textit{canonical embedding} in \cite{25}. In this paper, the author decides to retain this terminology, but also emphasizes that a "canonical embedding" is not an embedding in the literal sense of the definition in Section 2. In fact, it is virtually a planar diagram representing an embedding $\Gamma: G\mapsto \mathbb{RP}^2$.

\begin{figure}[H]
\centering
\includegraphics[width=8cm]{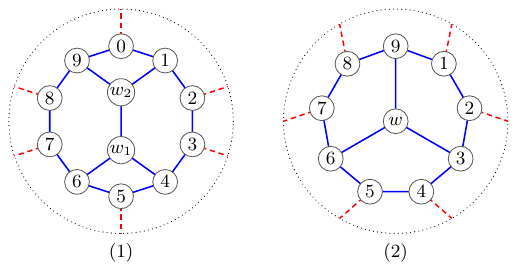}
\caption{signed graphs in $\mathcal{S}^*(3)$.}
\label{examples}
\end{figure}

As an example, it was proved in \cite{25} that $\mathcal{S}^*(3)$ consists of two members, and their canonical embeddings are shown in Figure \ref{examples}. In this paper, we always use blue solid and red dashed lines in figures to represent positive and negative edges respectively.
For $(G,\sigma)\in \mathcal{S}^*(k)$ with $E^-_\sigma(G)=\{x_1y_1,\ldots, x_ky_k\}$, we use $\mathcal{R}$ to denote the set $\{x_1,\ldots, x_k, y_1,\ldots, y_k\}$. The plane cubic graph in Theorem \ref{basis} is denoted by $G^\prime$. For a subgraph $H$ of $G$, the \textit{weight} of $H$ is defined as $\omega(H)=|V(H)\cap \mathcal{R}|$.

Throughout the remainder of the paper, all $(G, \sigma)\in \mathcal{S}^*(k)$ are deemed to be drawn with a canonical embedding. The signature $\sigma$ is thus minimum, i.e., the number of negative edges equals $k$.

It is well known that an edge cut $\partial_G(X)$ is minimal if and only if both $G[X]$ and $G[X^c]$ are connected. In $\mathbb{RP}^2$, for a minimal edge cut $\partial_G(X)$, there exists two compact subsurfaces $\Gamma_1$ and $\Gamma_2$ such that: (1) $\Gamma_1\cup\Gamma_2=\mathbb{RP}^2$; (2) $G[X]\subset \Gamma_1$, $G[X^c]\subset \Gamma_2$, and $G[X]\cap \Gamma_2=G[X^c]\cap \Gamma_1=\emptyset$; (3) $\Gamma_1\cap\Gamma_2$ is their common boundary, which is a single simple closed curve. Hence a minimal edge cut can be represented by a simple closed curve, and it can be encoded by a closed alternating sequence of consecutive edges and faces without repetition. For example, $B(\gamma)=[e_1F_1e_2\ldots e_rF_re_1]$ represents the closed curve $\gamma$ of a minimal cut formed by $r$ edges, where $e_i$ lies on the boundaries of its neighboring faces for $i\in [r]$. Conversely, an alternating sequence $[e_1F_1e_2\ldots e_rF_re_{r+1}]$ of consecutive edges and faces represents a curve or an arc, constructed by linking the midpoints of $e_i$ and $e_{i+1}$ with a simple arc within $F_i$ for every $i$.

For a signed graph $(G, \sigma)$ with a canonical embedding, we still use $e^-_i$ to denote the $k$ negative edges, and use $F^-_i$ to denote the $k$ faces lying in the crosscap in cyclic order ($i\in [k]$), where $e^-_i$ is in the common boundary of $F^-_i$ and $F^-_{i+1}$ for $i\in \textbf{Z}_k$. Deleting all $e^-_i$ and $F^-_j$ in $B(\gamma)$, the remaining subsequences correspond to a set of disjoint simple arcs $\alpha_1,\ldots,\alpha_s$ in the plane, which represent the positive edges contained in the edge cut represented by $\gamma$. We show that, in particular, if $\partial_G(X)$ is a minimal equilibrated cut, then there is exactly one such simple arc in its simple closed curve, that is, the simple closed curve of $\partial_G(X)$ traverses $G^\prime$ only once.
\begin{observation}\label{single}
Let $\partial_G(X)$ be a minimal equilibrated cut of $(G, \sigma)$ with simple closed curve $\gamma$ and 
$\{\alpha_1,\ldots,\alpha_s\}$ be the collection of disk arcs, then $s=1$.
\end{observation}
\begin{proof}
Suppose to the contrary that there are two or more simple arcs, then one of $G^\prime[X]$ and $G^\prime[X^c]$, say $G^\prime[X]$, is separated into $t(t \ge 2)$ connected components $ G^\prime[X_1], \cdots, G^\prime[X_t]$. These components are joined by only negative edges, hence $ d^-(X) < \sum_{i=1}^t d^-(X_i)$. On the other hand, $\sum_{i=1}^t d^+(X_i) = d^+(X)$ as there are no positive edges between different parts, and $d^+(X)=d^-(X)$ as $\partial_G(X)$ is an equilibrated cut. So $\sum_{i=1}^t d^+(X_i) < \sum_{i=1}^t d^-(X_i)$ and there exists some $i$ such that $d^+(X_i)<d^-(X_i)$, contrary to the premise that the signature $\sigma$ is minimum.
\end{proof}

Having built the relationship between edge cuts and simple closed curves, we prove the following property of canonical embedding of $(G,\sigma)\in \mathcal{S}^*(k)$. 
\begin{lemma}\label{key}
Let $\Gamma: G\mapsto \mathbb{RP}^2$ be the embedding of $(G,\sigma)\in \mathcal{S}^*(k)$ represented by its canonical diagram, and for $e\in E(G)$, let $\Gamma_e:G\setminus e\mapsto \mathbb{RP}^2$ be the restriction of $\Gamma$ to $G\setminus e$. The following statements hold.

(1) $\rho(\Gamma)=k$; (2) For every edge $e\in E(G)$, $\rho(\Gamma_e)<k$.
\end{lemma}
\begin{proof}
(1) It is obvious that there is an essential curve in the crosscap intersecting $\Gamma(G)$ at exactly the $k$ negative edges, so we only need to show $\rho(\Gamma)\geq k$.

Assume $\gamma$ is a minimal essential curve such that $|\gamma\cap \Gamma(G)|\leq k-1$. Since the fundamental group of $\mathbb{RP}^2$ is $\pi_1(\mathbb{RP}^2)\cong\textbf{Z}_2$, $\gamma$ intersects the crosscap border (the dotted circle in Figures) at an odd number of points. Because $G$ is cubic, a curve meeting a vertex can be perturbed locally to meet an incident edge instead, so we can modify $\gamma$ to be intersecting $\Gamma(G)$ at only internal points of edges in $E(\Gamma)$, and it admits an alternating edge-face sequence $B(\gamma)$. Moreover, the faces in $B(\gamma)$ are distinct, and thus so are edges: suppose a face $F$ is traversed twice by $\gamma$, then $B(\gamma)$ has the form $[\ldots e_1Fe_2\ldots e_3Fe_4\ldots]$. Take a point on each of the subcurves contained in $F$ and glue them into a single point $q$, then $\gamma$ is cut into two simple closed curves $\gamma_1$ and $\gamma_2$ by $q$, where $B(\gamma_1)=[\ldots e_1Fe_4\ldots]$ and $B(\gamma_2)=[\ldots e_2Fe_3\ldots]$. By the fact $\pi_1(\mathbb{RP}^2)\cong\textbf{Z}_2$, one of $\gamma_1$ and $\gamma_2$, say $\gamma_1$ is essential too, but $|\gamma_1\cap G|\leq |\gamma\cap G|-2$, contradicting the premise that $\gamma$ is minimal.

Let $\{\alpha_1,\ldots,\alpha_s\}$ be the simple arc set of $\gamma$ in the plane, and let $\{B_1,\ldots,B_s\}$ be their subsequences in $B(\gamma)$ respectively. Each $\alpha_i$ cuts the plane disk into two connected parts, hence corresponds to an edge cut $\partial_{G}(X_i)$ by constituting the positive edge set. We assert that there exists some $i\in[s]$ satisfying $|\alpha_i\cap\Gamma(G^\prime)|<\min \{\omega(X_i), \omega(X^c_i)\}$. If not, for each $i$, let $X_i$ be the part with smaller weight without loss of generality. In $B(\gamma)$, we replace every $B_i$ with the alternating edge-face sequence formed by negative edges in $\partial_G(X_i)\cap E_\sigma^-(G)$ together with the intervening crosscap faces. For example, assume the location of $B_1$ in $B(\gamma)$ is $[\ldots F^-_1B_1F^-_{t}\ldots]$, then $\partial_G(X_1)\cap E_\sigma^-(G)=\{e^-_1,\ldots,e^-_t\}$ or $\{e^-_{t+1},\ldots,e^-_k\}$, and after the replacement, $B(\gamma)$ becomes $[\ldots F^-_1e^-_1F^-_2\ldots e^-_{t-1}F^-_{t}\ldots]$ or $[\ldots F^-_1e^-_kF^-_k\ldots e^-_{t+1}F^-_{t}\ldots]$. This replacement preserves the intersection points of $\gamma$ and the crosscap border, as a result, we obtain a new edge-face sequence $B(\gamma^\prime)$ yielding a new essential curve $\gamma^\prime$ lying in the crosscap, see Figure \ref{facewidth}. $\gamma^\prime$ intersects $\Gamma(G)$ at only negative edges, and at no more than $k-1$ points by our assumption, which is impossible for such a curve. So there exists some $i\in[s]$ such that $|\alpha_i\cap\Gamma(G^\prime)|<\omega(X_i)$, and this implies $d_{(G,\sigma)}^+(X_i)<d_{(G,\sigma)}^-(X_i)$. By switching at $X_i$, we obtain another signature $\sigma^\prime$ with fewer than $k$ negative edges, which is a contradiction to the minimality of $\sigma$.
\begin{figure}[H]
\centering
\includegraphics[width=5cm]{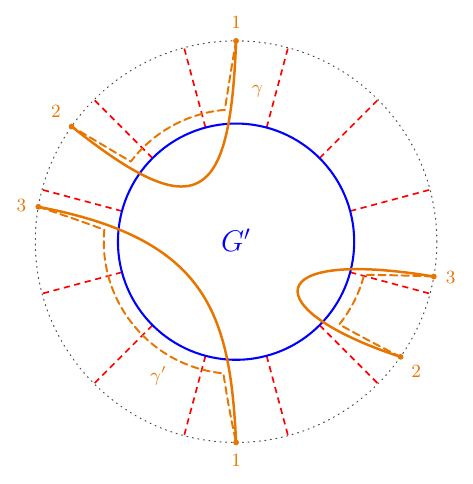}
\caption{}
\label{facewidth}
\end{figure}
(2) From the arguments in (1), it suffices to show that for every edge $e$ in $E(G)$, there exists a minimum signature $\pi$ switching equivalent to $\sigma$ such that $e\in E^-_{\pi}(G)$, and $(G, \pi)$ admits a canonical embedding under $\Gamma$.
\begin{observation}\label{matching}
For any equilibrated cut $\partial_G(X)$, the positive edges in $\partial_G(X)$ form a matching.   
\end{observation}
\begin{proof}
Assume to the contrary that an equilibrated cut $\partial_G(X)$ contains two positive edges incident to a vertex $w$. Because $G$ is cubic, after switching at $X$ and then switching at $w$, the resulting signature contains at most $k-1$ negative edges, contrary to the premise that $\sigma$ is minimum.
\end{proof}
We adopt a planar view of $(G,\sigma)$ based on $G^\prime$, in which every negative edge $x_iy_i$ in $E^-_\sigma(G)=\{x_1y_1,\ldots, x_ky_k\}$ is regarded to be broken by the crosscap into two halves on the plane, one incident to $x_i$ and the other to $y_i$. Let $\partial_G(Y)$ be a minimal equilibrated cut containing $e$. By Observation \ref{single}, the simple closed curve of $\partial_G(Y)$ in $\mathbb{RP}^2$ yields a single simple arc $\alpha$ in the plane. We switch at $Y$ and obtain another minimum signature $\pi$. On the plane, this operation is shown visually by the following procedure:
\begin{enumerate}

\item Cut along $\alpha$ to separate $G$ into two pieces induced by $G[Y]$ and $G[Y^c]$. Every positive edge in $\partial_G(Y)$ is regarded to be broken into two halves by the cut.

\item Invert the piece $G[Y]$ and move it to the opposite side of $G[Y^c]$.

\item Rejoin every pair of halves belonging to $\partial_G(Y)\cap E_\sigma^-(G)$ into a complete edge.

\item Change the sign of every edge in $\partial_G(Y)$.

\end{enumerate}
\begin{figure}[H]
\centering
\includegraphics[width=12.7cm]{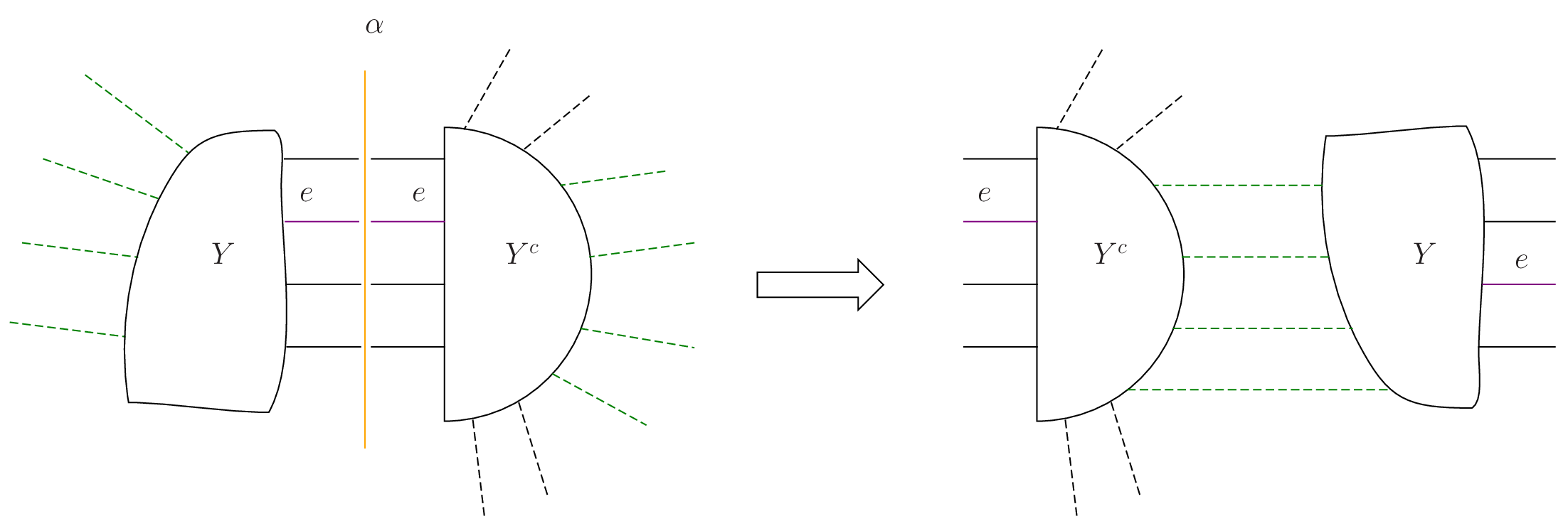}
\caption{cut-paste procedure.}
\label{planecut}
\end{figure}
The process is displayed in Figure \ref{planecut}, where we use green to highlight the negative edges belonging to $\partial_G(Y)$, and use purple to highlight $e$. All other edges are drawn in black in this figure to avoid making it too colorful. By Observation \ref{matching}, we finally obtain a drawing of canonical embedding of $(G,\pi)$, without changing $\Gamma: G\mapsto \mathbb{RP}^2$. This completes the proof of (2).
\end{proof}
An embedding $\Gamma: G\mapsto \Sigma$ is a \textit{minor minimal $k$-representative embedding} if $\rho(\Gamma)\geq k$ and for every proper minor $H$ of $G$, the representativity of the induced embedding of $H$ is less than $k$. The minor minimal $k$-representative embeddings in $\mathbb{RP}^2$ were characterized by S. Randby in 1997 \cite{key}.

\begin{figure}[H]
\centering
\includegraphics[width=5cm]{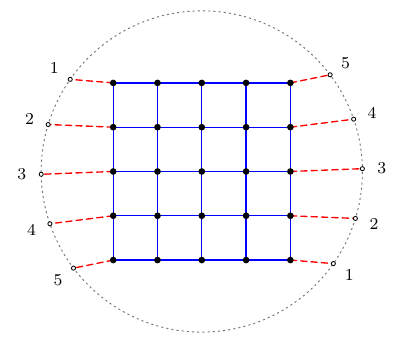}
\caption{$5\times 5$ M\"{o}bius grid.}
\label{MG_5}
\end{figure}

A \textit{$k\times k$ M\"{o}bius grid} is a $k\times k$ grid plus $k$ extra edges matching vertices in the first and last columns reversely. We use $MG_k$ to denote it. The $k$ twisted edges are assigned negative signs, and thus the $k\times k$ M\"{o}bius grid admits a "canonical embedding" drawing in $\mathbb{RP}^2$, Figure \ref{MG_5} shows an example of $MG_5$. We denote by $\Phi_k$ the embedding represented by this diagram: $\Phi_k: MG_k\mapsto\mathbb{RP}^2$, and let $\mathcal{T}_k$ be the delta-wye equivalence class of $\Phi_k$.
\begin{theorem}\label{randby}
(Randby\cite{key}) In $\mathbb{RP}^2$, an embedding $\Gamma$ is a minor minimal $k$-representative embedding if and only if it is in $\mathcal{T}_k$. 
\end{theorem}

Theorem \ref{randby} provides a characterization at the level of embeddings. We next show that, for $k\geq 4$, every graph underlying an embedding in $\mathcal T_k$ is uniquely embedded in $\mathbb{RP}^{2}$. This allows us to regard $\mathcal T_k$ as a class of graphs and subsequently discuss projective signatures relative to their uniquely determined embeddings. We begin with the celebrated Robertson–Vitray uniqueness theorem.

\begin{lemma}\label{unique}
(Robertson–Vitray\cite{vitray}) Let $\Sigma$ be a non-orientable surface and $\Gamma: G\mapsto \Sigma$ be an embedding of a $3$-connected graph $G$. If $\rho(\Gamma)\geq \widetilde{\epsilon}(\Sigma)+3$ when $\Sigma$ is nonorientable, then $\Gamma$ is the unique embedding of $G$ in $\Sigma$.
\end{lemma}

\begin{proposition}
For $k\geq 4$, graphs in $\mathcal{T}_k$ are uniquely embedded in $\mathbb{RP}^2$.
\end{proposition}
\begin{proof}
Since $MG_k$ is $3$-connected, to apply Lemma \ref{unique}, it suffices to show in $\mathcal{T}_k$, $3$-connectivity is preserved under a single $\Delta Y$ or $Y\Delta$ operation.

Let $G_{\triangle}\in\mathcal{T}_k$ contain a triangle $x-y-z-x$, and let $G_Y\in\mathcal{T}_k$ be obtained by replacing the triangle with a triad centered at a new vertex $v$. Suppose first that $G_Y$ is $3$-connected. For any $X\subseteq V(G_{\triangle})$ with $|X|\leq 2$, a path in $G_Y\setminus X$ passing through the triad can be converted into a walk in $G_\triangle\setminus X$ passing through a single edge in the triangle, for example, the subpath $x-v-y$ can be replaced by the edge $xy$. Hence $G_\triangle\setminus X$ is connected. Conversely, suppose $G_{\triangle}$ is $3$-connected. Let $X$ be a vertex set in $G_X$ with $|X|\leq 2$. If $v\notin X$, inversely, a path in $G_\triangle\setminus X$ passing through the triangle can be converted into a walk in $G_Y\setminus X$. If $v\in X$, let $X=\{v, w\}$, where $w$ is allowed to be absent. We assert $G_\triangle\setminus\{xy,yz,xz,w\}$ is connected, to which the connectivity of $G_Y\setminus X$ is equivalent. If not, then one of the vertices in $\{x, y, z\}$, say $x$, is separated from the other two by cut vertex $w$ of $G_\triangle\setminus\{xy,yz,xz\}$. Since $G_Y$ admits a minor-minimal $k$-representative embedding by Theorem \ref{randby}, $deg_{G_Y}(x)\geq 3$, hence $deg_{G_\triangle}(x)\geq 4$ and the part containing $x$ has a third vertex other than $\{x,w\}$, which implies $G_\triangle\setminus\{x, w\}$ is not connected, a contradiction to the premise. So $3$-connectivity propagates along $\Delta Y$ and $Y\Delta$ operations, and all graphs in $\mathcal{T}_k$ are uniquely embedded in $\mathbb{RP}^2$ for $k\geq 4$.
\end{proof}

\begin{figure}[H]
\centering
\includegraphics[width=12cm]{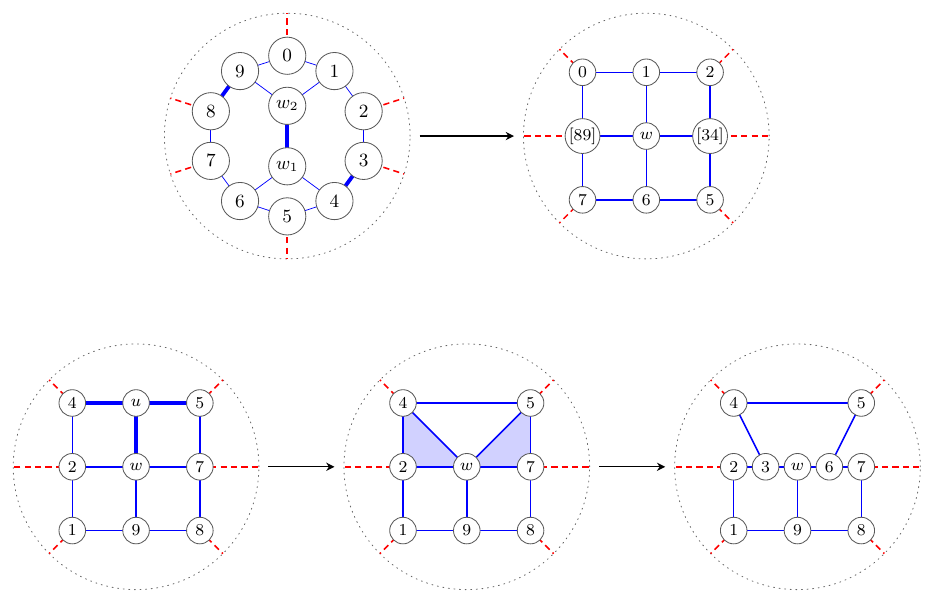}
\caption{Examples in $\mathcal{S}^*(3)$.}
\label{S3_examples}
\end{figure}

\begin{figure}[H]
\centering
\includegraphics[width=8cm]{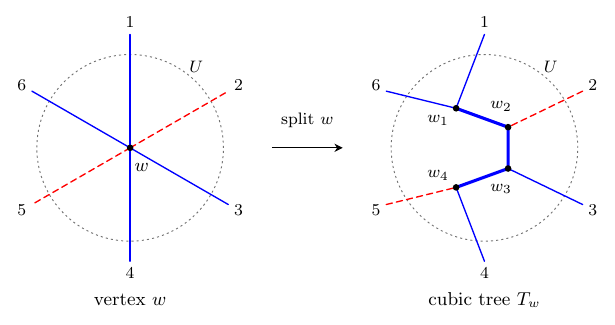}
\caption{tree expansion.}
\label{tree_expansion}
\end{figure}

From now on, for every graph $H\in \mathcal{T}_k$($k\geq 4$), we use $\Psi(H)$ to denote the unique embedding of $H$ in $\mathbb{RP}^2$. Using the canonical embedding drawing of $MG_k$ like the example in Figure \ref{MG_5}, one can construct a $k$-frustrated signed graph $(MG_k, \pi)$ by labeling exactly the $k$ edges in the crosscap negative, hence a cycle is negative precisely when its embedded image is essential, from the fact that $\pi_1(\mathbb{RP}^2)\cong\textbf{Z}_2$. In Subsection 3.2, we prove that the two descriptions are actually equivalent for a $k$-frustrated signed graph $(H,\pi)$ of $\mathcal{T}_k$, with the unique embedding $\Psi(H)$.
\begin{definition}
Fix an embedding $\Gamma: H\mapsto \mathbb{RP}^2$, and let $\pi$ be a signature of $H$. If in $(H,\pi)$, a cycle is unbalanced if and only if its image of $\Gamma$ is an essential cycle, then we call $\pi$ a \textit{projective signature}, and $(H,\pi)$ a projective signed graph.
\end{definition}
\begin{proposition}\label{prosign}
Let $k\geq 4$, and $(H,\pi)$ be a $k$-frustrated signed graph of $H\in \mathcal{T}_k$. Then $\pi$ is a projective signature if and only if every minimum signature $\pi'$ switching equivalent to $\pi$ satisfies the following conditions: (1) $E^-_{\pi'}(H)$ is a matching; (2) $(H,\pi')$ has a canonical embedding such that the $k$ edges in the crosscap are exactly $E^-_{\pi'}(H)$, and the positive part $H\setminus E^-_{\pi'}(H)$ is delta-wye equivalent to the $k\times k$ plane grid.
\end{proposition}
By Lemma \ref{key}, one can immediately deduce that for $(G,\sigma)\in \mathcal{S}^*(k)$, there exists an edge subset $E'\subseteq E(G)$ such that $G/E'$ is a graph with an embedding in $\mathcal{T}_k$. With Theorem \ref{essence1}, we have managed to make this process more specific, by showing that it only happens within the positive part of $(G,\sigma)$, and $E'$ is acyclic, hence the sets of positive and negative cycles of $(G,\sigma)$ are not changed by the contraction. As an example, Figure \ref{S3_examples} shows the contractions from the two signed graphs in $\mathcal{S}^*(3)$ to $\mathcal{T}_3$. The first signed graph is changed into $MG_3$ by contracting a matching $\{34, 89, w_1w_2\}$; the second signed graph is exactly the Petersen graph and already belongs to $\mathcal{T}_3$, thus $E'=\emptyset$ and we display the $\Delta Y$ and $Y\Delta$ operations transforming $MG_3$ into it.
Conversely, for an embedding $\Gamma: (G,\sigma)\mapsto \Sigma$, a \textit{positive tree expansion} of $\Gamma$ is an embedded cubic signed graph $(\tilde{G},\sigma)$ obtained by splitting every vertex $w$ of degree at least $4$ into an all-positive cubic tree $T_w$ within a local disk, while the signature and local cyclic order of original edges in $E_G(w)$ are not changed, see Figure \ref{tree_expansion}. So every $(G,\sigma)\in \mathcal{S}^*(k)$ is a positive tree expansion of a $k$-frustrated projective signed graph of $\mathcal{T}_k$. Theorem \ref{essence2} shows that graphs in $\mathcal{S}^*(k)$ and $\mathcal{T}_k$ are uniquely embedded in $\mathbb{RP}^2$, Thus the definition can be made independently of the choice of embedding.
The proofs of Theorem \ref{essence1} and Theorem \ref{essence2} are displayed in the following two subsections.

\subsection{$\Delta$Y/Y$\Delta$ operations within the disk}
We will use the medial graph theory developed by E.B.Curtis et al. \cite{medial} to show that the edge contractions and $\Delta$Y/Y$\Delta$ operations in Theorem \ref{essence1} are performed in the disk part of $\mathbb{RP}^2$, with no negative edge of $(G, \sigma)$ involved.

A \textit{circular plane graph} $P$ is a graph embedded in a disk $D$ with boundary circle $\mathcal{C}$, so that $P \cap \mathcal{C}$ is exactly a vertex subset of $V(P)$, called \textit{nodes}, and the rest of $P$ is in the interior of $D$. Suppose $P$ has $2k$ nodes $\textbf{n}_1,\ldots,\textbf{n}_{2k}$ arranged in clockwise order. For each $e\in E(P)$, let $m_e$ be its midpoint, the \textit{medial graph} of $P$, denoted by $\mathscr{M}(P)$, is defined as follows:

(1) Place $4k$ points $\textbf{t}_1,\ldots,\textbf{t}_{4k}$ on $\mathcal{C}$ so that
\begin{displaymath}
\textbf{t}_1<\textbf{n}_1<\textbf{t}_2<\textbf{t}_3<\textbf{n}_2\cdots <\textbf{t}_{4k-1}<\textbf{n}_{2k}<\textbf{t}_{4k}<\textbf{t}_1
\end{displaymath}
in the clockwise circular order around $\mathcal{C}$;

(2) Let $V(\mathscr{M}(P))=\{m_e\}_{e\in E(P)} \cup \{\textbf{t}_i\}_{i=1}^{4k}$. For two edges $e_1, e_2\in E(P)$, $m_{e_1}$ and $m_{e_2}$ are joined by an edge whenever $e_1$ and $e_2$ have a common vertex and incident to the same face in the embedding of $P$;

(3) For $j\in [2k]$, the point $\textbf{t}_{2j}$ is joined by an edge to $m_e$, where $e$ is the edge coming first after boundary interval $\textbf{n}_j\textbf{t}_{2j}$ in clockwise order around $\textbf{n}_j$. The point $\textbf{t}_{2j-1}$ is joined by an 
edge to $m_{e'}$, where $e'$ is the edge coming first after interval $\textbf{n}_j\textbf{t}_{2j}$ in counterclockwise order around $\textbf{n}_j$.

\begin{figure}[H]
\centering
\includegraphics[width=12cm]{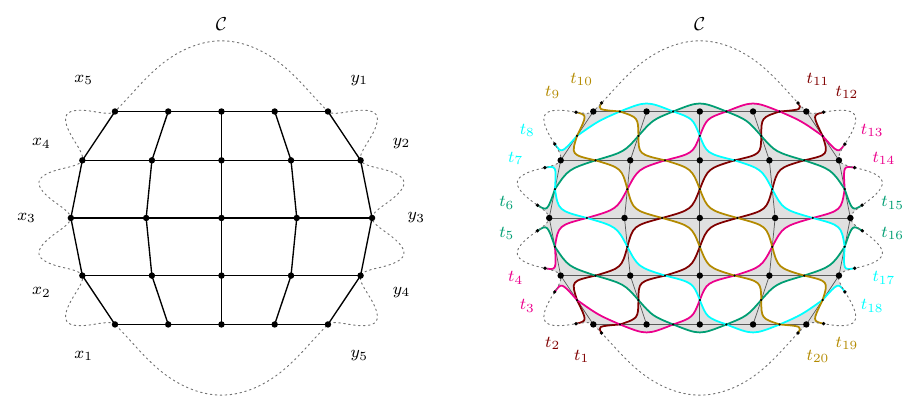}
\caption{circular plane graph $GR_5$ and $\mathscr{M}(GR_5)$.}
\label{medial}
\end{figure}

In $\mathscr{M}(P)$, from the definition, every $m_e$ is a vertex of degree four, and every $\textbf{t}_j$ is a pendant vertex. A \textit{geodesic} is a maximal walk that passes through every $m_e$ in such a way that, it leaves along the opposite incident edge upon entering $m_e$. Thus a geodesic is a closed walk containing only midpoints $m_e$, or has two ends $\textbf{t}_{j_1}$ and $\textbf{t}_{j_2}$. Figure \ref{medial} shows the medial graph of $5\times 5$ plane grid, where the nodes are $\{x_1,\ldots, x_5, y_1,\ldots, y_5\}$, and $\mathcal{C}$ is constructed by adding external simple arcs between $\{x_i, x_{i+1}\}$, $\{y_i, y_{i+1}\}$, $\{x_1, y_5\}$ and $\{x_5, y_1\}$ ($1\leq i\leq 4$).
Edges and geodesics in $\mathscr{M}(P)$ are regarded to be topological curves in the disk $D$, and $D$ is shattered into two different types of regions, one containing exactly one vertex in $V(P)$, and one corresponding to a face of $P$. So the regions of $D\setminus \mathscr{M}(P)$ have a natural coloring: regions containing vertices of $P$ are colored black, and regions corresponding to faces of $P$ are colored white. Every medial edge separates a black region from a white region. The boundary circle $\mathcal{C}$ is likewise shattered into two types of segments: for $i=1,\ldots,2k \mod 4k$, interval $(\textbf{t}_{2i-1}, \textbf{t}_{2i})$ is on the boundary of a black region, and $(\textbf{t}_{2i}, \textbf{t}_{2i+1})$ is on the boundary of a white region. We call the embedded medial graph together with this distinguished coloring the \textit{chessboard} of $P$.

Conversely, the circular plane graph $P$ can be reconstructed from its chessboard as follows. Place one vertex $v_\textbf{B}$ in every black region $\textbf{B}$, taking $v_\textbf{B}=\textbf{n}_i$ when $\textbf{B}$ contains the designated boundary interval $(\textbf{t}_{2i-1},\textbf{t}_{2i})$. Join an edge through two vertices $v_{\textbf{B}_1}$ and $v_{\textbf{B}_2}$ whenever the two black regions $\textbf{B}_1$ and $\textbf{B}_2$ are opposite. Thus these two constructions between $P$ and $\mathscr{M}(P)$ are mutually inverse up to ambient isotopy of $D$.

\begin{figure}[H]
\centering
\includegraphics[width=12cm]{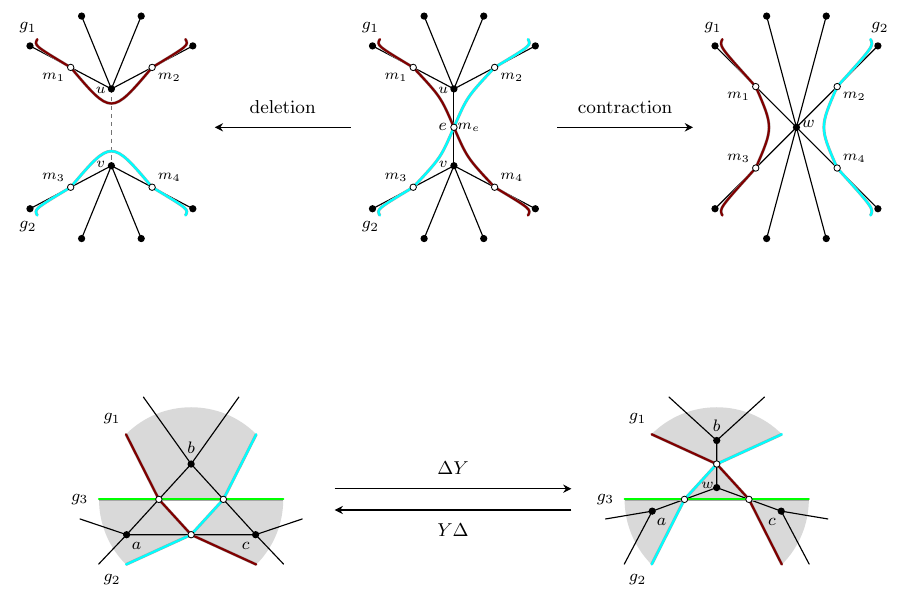}
\caption{geodesic modifications.}
\label{motion}
\end{figure}

In particular, edge removal (deletion or contraction) and $\Delta$Y/Y$\Delta$ operations in $P$ can be expressed in $\mathscr{M}(P)$ by modifications of geodesics. Removing edge $e\in E(P)$ corresponds to breaking the crossing point $m_e$ in $\mathscr{M}(P)$ and then smoothing the four branches. Suppose three geodesics cross pairwise and the subarcs between their three crossings bound an empty triangular region. A \textit{motion} is to move one of the three geodesics across the empty triangular region within a local disk to obtain another pairwise crossing configuration bounding an empty triangular region of opposite color, while the order of the two crossings in each geodesic is reversed, and all geodesics remain unchanged outside the local disk. These geodesic modifications are illustrated in Figure \ref{motion}.
A medial graph is said to be \textit{lensless} if every geodesic begins and ends at distinct boundary points $\textbf{t}_{j_1}$ and $\textbf{t}_{j_2}$, has no self-intersection, and any two geodesics have at most one intersection point. If a medial graph $\mathscr{M}(P)$ is lensless with geodesics $\{g_1,\ldots,g_{2k}\}$, we relabel $\{\textbf{t}_1,\ldots,\textbf{t}_{4k}\}$ with the geodesic indices and obtain an integer sequence $\mathscr{A}(P)=z_1,\ldots,z_{4k}$, in which each of the numbers in $[2k]$ occurs exactly twice, and two geodesics $g_i$ and $g_j$ intersects if and only if the two occurrences of $i$ alternate with the two occurrences of $j$ in $\mathscr{A}(P)$. The set of crossing pairs is therefore determined by $\mathscr{A}(P)$, and the the remaining ambiguity in reconstructing $P$ is that the positions of these points within $D$ may differ. Curtis et al. \cite{medial} proved that given a cyclic sequence $\mathscr{A}(P)$, every lensless medial graph with that endpoint pairing can be transformed into a unique standard medial graph by only geodesic motions. We refer the reader to \cite{medial} for the details.

Because geodesic motion corresponds to $\Delta Y$ or $Y\Delta$ operation in the original circular plane graph, and delta-wye equivalence is an equivalence relation, the following lemma holds.

\begin{lemma}\label{lenseq}
Suppose $P_1$ and $P_2$ are two connected circular plane graphs, each with the same labeled boundary nodes in the same cyclic order. Assume both $\mathscr{M}(P_1)$ and $\mathscr{M}(P_2)$ are lensless. Then $P_1$ and $P_2$ are delta-wye equivalent if and only if $\mathscr{A}(P_1)=\mathscr{A}(P_2)$.
\end{lemma}

Now we proceed to the main part of the proof of Theorem \ref{essence1}. Let $GR_k$ be the $k\times k$ plane grid, where the vertices in the first column are labeled with $x_1,\ldots,x_k$, and the vertices in the last column are labeled with $y_1,\ldots,y_k$ in reverse order. Embed $GR_k$ in the disk $D$ and let $\{x_1,\ldots,x_k, y_1,\ldots,y_k\}$ be nodes on boundary $\mathcal{C}$. Insert $4k$ points $\textbf{t}_1,\ldots,\textbf{t}_{4k}$ on $\mathcal{C}$ so that $\textbf{t}_{2i-1}< x_i<\textbf{t}_{2i}$ and $\textbf{t}_{2k+2i-1}< y_i<\textbf{t}_{2k+2i}$ for $i\in [k]$.

\begin{fact}
$\mathscr{M}(GR_k)$ is lensless, and the ends of geodesics are $\{\textbf{t}_{2i-1}, \textbf{t}_{2k+2i}\}$ and $\{\textbf{t}_{2i}, \textbf{t}_{2k+2i-1}\}$ ($i\in [k]$).
\begin{proof}
It can be checked directly from the definition that pairs $\{\textbf{t}_{1}, \textbf{t}_{2k+2}\}$ and $\{\textbf{t}_2, \textbf{t}_{2k+1}\}$ are linked by two self-intersection free geodesics $g_1$ and $g_2$ respectively, and $g_1$ and $g_2$ are disjoint. By symmetry of $GR_k$, all other geodesics are obtained by translating $g_1$ and $g_2$, so all geodesics are simple, joining the pairs $\{\textbf{t}_{2i-1}, \textbf{t}_{2k+2i}\}$ and $\{\textbf{t}_{2i}, \textbf{t}_{2k+2i-1}\}$ for $i\in [k]$, and any two geodesics intersect at no more than one point. 
\end{proof}
\end{fact}

Let $Q$ be a circular plane graph with the same nodes $\{x_1,\ldots,x_k, y_1,\ldots,y_k\}$. Adding $k$ negative edges $\{x_1y_1,\ldots,x_ky_k\}$ reversely in the crosscap, we then obtain a projective plane signed graph $\tilde{Q}$ with an embedding $\Gamma$. We say $\tilde{Q}$ satisfies positive-minor criticality, denoted by $\mathcal{P}$, if $\rho(\Gamma(\tilde{Q}))=k$, and the representativity of every proper minor retaining the $k$ negative edges is strictly less than $k$. In the following, we will show that $Q$ is delta-wye equivalent to $GR_k$ if $\mathcal{P}$ holds for the signed graph $\tilde{Q}$.

\begin{proposition}\label{q1}
$\mathscr{M}(Q)$ is lensless.
\end{proposition}
\begin{proof}
We first point out that $\mathcal{P}$ propagates along $\Delta Y$ and $Y\Delta$ operations within the ambient disk $D$ containing $Q$, that is to say, no negative edge is involved.

It suffices to show $\mathcal{P}$ holds after a single $\Delta Y$ or $Y\Delta$ operation. Let $Y(\tilde{Q})$ be the signed graph obtained by applying one $\Delta Y$ operation on $Q$, and $U\subseteq D$ be the local disk in the definition. Arbitrarily select a positive edge $e$ of $Y(\tilde{Q})$. If $e$ is not in $U$, then $e$ is also an edge of $Q$. Since every minimal essential curve intersects the triangle or the triad at the same number of points when it passes through $U$, by $\mathcal{P}$ of $\tilde{Q}$, $\rho(Y(\tilde{Q}))=k$, and either deleting or contracting $e$ will decrease the representativity. If $e$ lies in $U$, then $e$ is an edge of the triad. Denote the three vertices of the triangle by $x$, $y$ and $z$, and the center of the triad by $v$. Choose $e=xv$ without loss of generality, then $Y(\tilde{Q})\setminus e$ is isomorphic to the graph obtained from $\tilde{Q}$ by deleting $xy$, $xz$, and then subdividing the edge $yz$ with a new vertex $v$, while $Y(\tilde{Q})/e$ is isomorphic to $\tilde{Q}-yz$. From $\mathcal{P}$ of $\tilde{Q}$, representativity is strictly decreased in both cases. Let $\triangle(\tilde{Q})$ be the signed graph obtained by a single $Y\Delta$ operation on $Q$. For the same reason, $\rho(\triangle(\tilde{Q}))=k$ and deleting or contracting any positive edge will decrease the representativity. The verification is almost identical and we omit the details here. 

From now on we verify the lenslessness of $\mathscr{M}(Q)$. First, assume there is a self-intersecting or closed geodesic $g$, then in both cases a domain homeomorphic to the disk is yielded, bounded by a closed curve lying in $g$. From all such domains, we choose an inclusion-minimal one, and denote it by $\Omega$, then the boundary of $\Omega$ is a simple closed curve. $\Omega$ is the union of several black and white regions, we claim that the domain $\Omega$ can be transformed into a single region via a series of geodesic motions, called \textit{evacuation} of $\Omega$, executed in the following process:

\textbf{Evacuation}: if there exist geodesics $h_1$, $h_2$ intersecting at a point in the interior of $\Omega$, the point induces two sectors formed by $h_1$, $h_2$ and $g$. From such geodesic pairs $(h_1, h_2)$, choose the one that: (1) the intersection point on $h_1$ is the closest one to $g$; (2) subject to (1), one of the two sectors contains the fewest regions. Then expel the intersection point from $\Omega$ with a motion of $\{h_1, h_2, g\}$. Repeat the above steps until no geodesics intersect within $\Omega$. This procedure is illustrated in Figure \ref{evacuation}, in which we shade in green the potential motion sectors in the initial stage.
\begin{figure}[H]
\centering
\includegraphics[width=12cm]{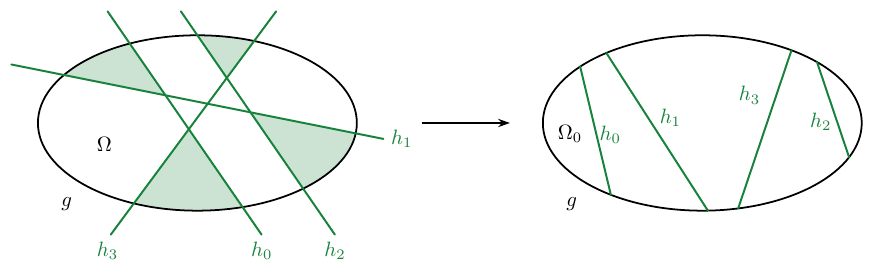}
\caption{evacuation.}
\label{evacuation}
\end{figure}
The domain $\Omega$ is empty after \textbf{evacuation}. If not, the geodesic segments contained in $\Omega$ are disjoint, from which we take an outermost one $h_0$, hence $h_0$ and $g$ form the boundary circle of a single region $\Omega_0$. Because two geodesics cannot be tangent at their intersection point by definition, there is no intermediate vertex in $V(\mathscr{M}(Q))$ between the two fixed intersection points of $h_0$ and $g$, that is to say, the boundary circle of $\Omega_0$ contains exactly two vertices of $\mathscr{M}(Q)$. By the correspondence between $Q$ and its chessboard, \textbf{evacuation} is interpreted in $Q$ in the following way: after a series of $\Delta Y$ and $Y\Delta$ operations within $D$, a subdivision vertex ($\Omega_0$ is black) or a pair of parallel edges ($\Omega_0$ is white) appears, either of which is a contradiction to $\mathcal{P}$. So $\Omega$ is a single region. Similarly, the boundary circle of $\Omega$ contains exactly one vertex, then $\Omega$ corresponds to a pendant vertex ($\Omega$ is black) or a loop ($\Omega$ is white), which is impossible again by $\mathcal{P}$.

Now all geodesics in $\mathscr{M}(Q)$ are simple arcs joining a pair $\{t_i, t_j\}$ on $\mathcal{C}$, we then prove that two geodesics will not intersect at two points. On the contrary, let $(g_1, g_2)$ be such a pair, that the domain $\Omega'$ enclosed by them contains the minimum number of regions in the chessboard, then any third geodesic $h$ passing through $\Omega'$ intersects $g_1$ and $g_2$ at exactly one point respectively. In this case, \textbf{evacuation} can be executed too on $\Omega'$ by modifying $g$ to one of $g_1$ and $g_2$ in the process, and the geodesic segments contained in $\Omega'$ are arranged sequentially between the two intersection points of $g_1$ and $g_2$ after \textbf{evacuation}, say $\{h'_1,\ldots,h'_s\}$. We then expel these segments from $\Omega'$ consecutively via geodesic motions on $\{g_1, g_2, h'_i\}$. In the end, $\Omega'$ is transformed into a single region, which is a contradiction from the above explanation. This completes the proof of lenslessness of $\mathscr{M}(Q)$.
\end{proof}

\begin{proposition}\label{q2}
$\mathscr{A}(Q)=\mathscr{A}(GR_k)$. That is, the ends of geodesics in $\mathscr{M}(Q)$ are $\{\textbf{t}_{2i-1}, \textbf{t}_{2k+2i}\}$ and $\{\textbf{t}_{2i}, \textbf{t}_{2k+2i-1}\}$ ($i\in [k]$).
\end{proposition}
\begin{proof}
Fix an arbitrary $j$ in $1,\ldots,2k \mod 4k$, place two points $\textbf{a}$, $\textbf{b}$ on $\mathcal{C}$ such that $\textbf{t}_{2j}<\textbf{a}<\textbf{t}_{2j+1}$ and $\textbf{t}_{2k+2j}<\textbf{b}<\textbf{t}_{2k+2j+1}$. Both points $\textbf{a}$ and $\textbf{b}$ lie on the boundary of a white region, and divide $\mathcal{C}$ into two halves, each of which contains exactly $k$ consecutive nodes as well as the $2k$ $\textbf{t}_i$'s around them. Let $\alpha$ be a simple arc linking $\textbf{a}$ and $\textbf{b}$ such that $\alpha$ passes through no vertex in $V(\mathscr{M}(Q))$ and $|\alpha\cap \mathscr{M}(Q)|$ is minimized.
\begin{claim}\label{uninter}
$\alpha$ intersects every geodesic at most once.
\end{claim}
\begin{proof}
On the contrary, suppose $\alpha$ intersects a geodesic $g$ at two distinct points $p_1$ and $p_2$. Let $\alpha(p_1, p_2)$ and $g(p_1, p_2)$ be the segments of $\alpha$ and $g$ between $p_1$ and $p_2$ respectively, then a domain $\Omega$ homeomorphic to a disk is bounded by them. Choose another point $p_1'\,\,(p_2')\notin V(\mathscr{M}(Q))$ on $\alpha\setminus \Omega$ arbitrarily close to $p_1$ ($p_2$), then we construct a new simple curve $\alpha'$ by replacing $\alpha(p_1, p_2)$ with a curve very close to $g$ linking $p_1'$ and $p_2'$, say $\alpha'(p_1', p_2')$, see Figure \ref{probe}.
\begin{figure}[H]
\centering
\includegraphics[width=8cm]{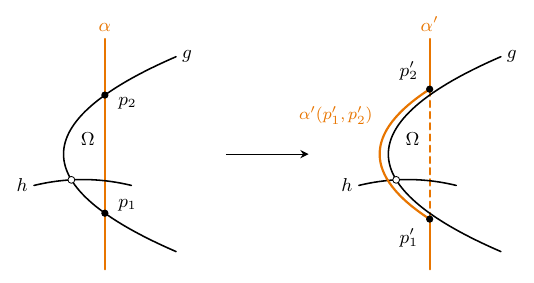}
\caption{detour of the probe.}
\label{probe}
\end{figure}
Let us analyze the intersections of $\alpha'$ and geodesics. Compared to $\alpha$, only intersections along the replacement arc need to be considered. Let $h$ be a geodesic intersecting $\alpha'$ at a point in $\alpha'(p_1', p_2')$. Since $\alpha'(p_1', p_2')$ is close enough to $g(p_1, p_2)$, $h$ must pass through $\Omega$ via a point in $g$, and lenslessness gives $|h\cap \alpha'(p_1', p_2')|=|h\cap g|=1$. Therefore $h$ must cross $\alpha(p_1, p_2)$ when it leaves $\Omega$, which implies $|h\cap \alpha'(p_1', p_2')|\leq |h\cap \alpha(p_1, p_2)|$. Consequently, the replacement eliminates at least two intersection points $p_1$ and $p_2$ with geodesics in $\mathscr{M}(Q)$, and this contradicts the choice of $\alpha$.
\end{proof}

We will show that the ends of every geodesic must lie in different halves by using $\alpha$ as a probe. Since $\alpha$ traverses no vertex in $\mathscr{M}(Q)$, its trail in the chessboard must be
\begin{center}
   white $\rightarrow$ black $\rightarrow$ white $\rightarrow$ $\cdots$ black $\rightarrow$ white,
\end{center}
and each inner region in the above sequence is intersected exactly twice on the boundary. In terms of the signed graph $\tilde{Q}$, the simple arc $\alpha$ traversing disk $D$ passes through the plane part $Q$ via only vertices in $V(Q)$, and it can be extended in the crosscap to an essential curve of $\mathbb{RP}^2$ without introducing additional intersections. Since $\rho(\tilde{Q})=k$, $|\alpha\cap V(Q)|\geq k$, equivalently, at least $k$ black regions of the chessboard are traversed in the trail of $\alpha$, and thus the total number of intersection points of $\alpha$ and geodesics is at least $2k$. By Claim \ref{uninter}, it is exactly $2k$, and $\alpha$ intersects all geodesics in $\mathscr{M}(Q)$. If one of the $2k$ geodesics has both ends in the same half of $\mathcal{C}$, then by planarity, $\alpha$ must be tangent to it at the unique intersection point, and that tangency can be removed by a local perturbation of $\alpha$. So each intersection point on $\alpha$ is a witness of a unique geodesic linking the two halves of $\mathcal{C}$ divided by $\textbf{a}$ and $\textbf{b}$, and the two ends of every geodesic in $\mathscr{M}(Q)$ are $\textbf{t}_{i_1}$ and $\textbf{t}_{i_2}$, where $2j+1\leq i_1\leq2j+2k$ and $2j+2k+1\leq i_2\leq 2j$ in the sense of modulo $4k$. To achieve this for an arbitrarily chosen $j$, the only solution is that the ends of every geodesic lie around a pair of antipodal nodes on $\mathcal{C}$, and thus geodesics in $\mathscr{M}(Q)$ are divided into $k$ pairs, each of which links the two sets $\{\textbf{t}_{2i-1}, \textbf{t}_{2i}\}$ and $\{\textbf{t}_{2k+2i-1}, \textbf{t}_{2k+2i}\}$ for some $i\in [k]$.

Within each pair, the ends of the two geodesics are $\{\textbf{t}_{2i-1}, \textbf{t}_{2k+2i}\}$ and $\{\textbf{t}_{2i}, \textbf{t}_{2k+2i-1}\}$, or $\{\textbf{t}_{2i-1}, \textbf{t}_{2k+2i-1}\}$ and $\{\textbf{t}_{2i}, \textbf{t}_{2k+2i}\}$. Call the pair is parallel for the former case and crossing for the latter. For every possible crossing pair, we delete their intersection point and then smooth the four branches, the ends of the geodesics are interchanged as a result and the pair is transformed into a parallel one. By planarity, geodesics from different pairs intersect at exactly one point, so there exist $4\binom{k}{2}$ such intersections, yielding $4\binom{k}{2}$ points in total. The transformation does not change the intersection relationship between geodesics from different pairs, nor the previous point set, so the $4\binom{k}{2}$ points witness all $4\binom{k}{2}$ geodesic intersections in the new medial graph, hence two geodesics intersect at no more than one point, and we finally obtain a lensless medial graph with the same geodesic endpoint pairing as $\mathscr{M}(GR_k)$. Let $Q_0$ be the circular plane graph reconstructed from it, and $\tilde{Q}_0$ be the signed graph by adding the $k$ negative edges to $Q_0$. By Lemma \ref{lenseq}, $Q_0$ and $GR_k$ are delta-wye equivalent, and thus $\tilde{Q}_0\in \mathcal{T}_k$ because no negative edge was involved in the $\Delta Y$ and $Y\Delta$ operations within the ambient disk $D$, so $\rho(\tilde{Q}_0)=\rho(MG_k)=\rho(\tilde{Q})=k$ by Theorem \ref{randby}. On the other hand, $\tilde{Q}_0$ is obtained from $\tilde{Q}$ by possibly deleting and contracting some positive edges from the correspondence between circular plane graph and its chessboard, but the property $\mathcal{P}$ of $\tilde{Q}$ rules out any such operation, so $Q_0=Q$, and all geodesic pairs in $\mathscr{M}(Q)$ are parallel, which is equivalent to $\mathscr{A}(Q)=\mathscr{A}(GR_k)$. 
\end{proof}

We can now apply the above propositions about $Q$ to obtain Proposition \ref{prosign} and Theorem \ref{essence1}.

\begin{proof}[\textbf{Proof of Proposition \ref{prosign}:}]
We only need to show necessity since sufficiency is obvious. Let $\pi'$ be a minimum signature switching equivalent to $\pi$. The boundary cycle of every face contains an even number of edges in $E^-_{\pi'}(H)$, so the dual edges corresponding to
the $k$ negative edges in $E^-_{\pi'}(H)$ form an Eulerian subgraph of $H^*$, the dual graph of $H$, and thus it has an edge-disjoint decomposition into dual cycles $C_1^*\cup\cdots C_s^*$. Every $C_i^*$ corresponds to an alternating sequence $B[C_i^*]$ of negative edges and faces in $H$ without repetition, so it can be represented by a simple closed curve $\gamma_i^*$ in $\mathbb{RP}^2$. The intersection points of distinct $\gamma_i^*$ and $\gamma_j^*$ are deemed to lie in the interior of faces of $H$ since $C_i^*$ and $C_j^*$ are edge disjoint.

We assert that $s=1$ and the unique cycle $C_1^*$ induces an essential curve $\gamma_1^*$ in $\mathbb{RP}^2$. Arbitrarily select an unbalanced cycle $C_0$ in $(H,\pi')$, and let $E^*_-(C_0)$ be the set of dual edges of $E(C_0)\cap E^-_{\pi'}(H)$. The curve $\gamma_i^*$ intersects $C_0$ if and only if $E(C_i^*)\cap E^*_-(C_0)\neq \emptyset$, and intersection points lie in the simple arcs representing negative edges of $C_0$. Then $|E(C_0)\cap E^-_{\pi'}(H)|=\sum_{i=1}^s|E(C_i^*)\cap E^*_-(C_0)|$ is odd because $C_0$ is unbalanced, so there exists an index $j$ such that $|E(C_j^*)\cap E^*_-(C_0)|$ is odd. The curve $\gamma_j^*$ must be essential, for otherwise $\gamma_j^*$ encloses a disk traversed by $C_0$, which implies $|E(C_j^*)\cap E^*_-(C_0)|$ is even because $C_0$ is essential. By the restriction $\rho(\Psi(H))=k$ from Theorem \ref{randby}, all $k$ dual edges corresponding to negative edges of $E^-_{\pi'}(H)$ all lie in $C_j^*$, and thus the decomposition consists of exactly one cycle $C_1^*$.

Moreover, if two edges of $E^-_{\pi'}(H)$, say $e^-_1$ and $e^-_2$ are incident to the same vertex $w$, we then replace the subsequence $B(C_1^*)=[\ldots e^-_1\ldots e^-_2\ldots]$ with $[\ldots w \ldots]$, visually, $\gamma_1^*$ is shortcut through vertex $w$ and a new essential curve intersecting $H$ at no more than $k-1$ points is yielded, which is a contradiction to Theorem \ref{randby}. Thus, we conclude that the negative edges of $E^-_{\pi'}(H)$ form a matching traversed by a minimal essential curve $\gamma_1^*$. Since $\gamma_1^*$ is isotopic to the core of a M\"{o}bius band in $\mathbb{RP}^2$, we can choose a disk-and-crosscap drawing of the projective plane, in which $\gamma_1^*$ intersects the crosscap border at exactly one point. As we have done in the proof of Lemma \ref{key}, by executing the cut-paste procedure shown in Figure \ref{planecut} with $\gamma_1^*$, we obtain a canonical embedding of $(H,\pi')$, with the $k$ negative edges of $E^-_{\pi'}(H)$ drawn through the crosscap in reverse matching order.

Finally, under the canonical embedding, a circular plane graph can be constructed from the positive part $H\setminus E^-_{\pi'}(H)$ by designating the end vertices of $E^-_{\pi'}(H)$ as boundary nodes. Take $Q=H\setminus E^-_{\pi'}(H)$ and $\tilde{Q}=(H,\pi')$ in Proposition \ref{q1} and \ref{q2}, since $(H,\pi')$ satisfies $\mathcal{P}$ by Theorem \ref{randby}, $H\setminus E^-_{\pi'}(H)$ is delta-wye equivalent to $GR_k$.
\end{proof}

\begin{proof}[\textbf{Proof of Theorem \ref{essence1}}]:
Given any canonical embedding of $(G, \sigma)\in \mathcal{S}^*(k)$, let $\Gamma: G\mapsto \mathbb{RP}^2$ be the graph embedding on $\mathbb{RP}^2$ reflected in it, and $G^+$ be the circular plane graph embedded in $D$ induced by $E_\sigma^+(G)$ with the same nodes $\{x_1,\ldots,x_k, y_1,\ldots,y_k\}$. From Lemma \ref{key}, it can be inferred directly by logic that there exists a subset $E'\subseteq E_\sigma^+(G)$ such that the induced embedding $\Gamma_{E'}: G/E'\mapsto \mathbb{RP}^2$ satisfies $\rho(\Gamma_{E'})=k$, and for any edge $e\notin E'$, $\rho(\Gamma_{E'\cup \{e\}})<k$. Let $Q=G^+/E'$ and $\tilde{Q}=G/E'$, then it is immediate that $G^+/E'$ is delta-wye equivalent to $GR_k$ and $G/E'$ is a $k$-frustrated projective signed graph of $\mathcal{T}_k$.

It remains to prove that $G[E']$ is acyclic. Suppose $G[E']$ contains a cycle $C$, then $C$ is a trivial curve in $\mathbb{RP}^2$, enclosing a disk $\Omega$. Arbitrarily take an edge $e$ on $C$. By Lemma \ref{key}, $e$ is traversed by a minimal essential curve $\gamma$ with $|\gamma\cap C|\geq 2$. As a result, $|\gamma\cap \Gamma_{E'}(\tilde{Q})|\leq k-1$, contrary to the premise that $\rho(\Gamma_{E'})=k$.
\end{proof}

\subsection{$\mathcal{S}^*(k)$ from tree expansion}
It is immediate from Theorem \ref{essence1} that every $(G,\sigma)\in \mathcal{S}^*(k)$ with $|E^-_\sigma(G)|=k$ is a positive tree expansion of a $k$-frustrated projective signed graph of $\mathcal{T}_k$. In this subsection, we will prove Theorem \ref{essence2} explaining this relationship in more depth, from which the two quadratic bounds of $|E(G)|$ in Corollary \ref{quad}, as well as the finiteness of $\mathcal{S}^*(k)$ follow.

\begin{proof}[\textbf{Proof of Theorem \ref{essence2}:}]

Let $T$ be a tree expansion of $G_0\in \mathcal{T}_k$. We decompose the vertex splittings into steps, one of which is the extension of a single edge:
\begin{displaymath}
G_0<\cdots<G_{t-1}<G_t=T,
\end{displaymath}

where $G_{i}$ is obtained from $G_{i-1}$ by splitting vertex $w_i$ into two vertices $\{u_i,v_i\}$, and adding a new edge $e_i=u_iv_i$ ($i=1,\ldots,t$). $G_0$ is non-planar and admit an embedding in $\mathbb{RP}^2$, so as every other $G_i$ in the chain. Moreover, by Lemma \ref{genus}, their embeddings are all 2-cell embeddings.

We assert that the uniqueness of embedding propagates along the chain. For $i=1,\ldots,t$, let $A_i$, $B_i$ be the fixed nonempty sets of neighboring vertices of $w_i$ assigned to $u_i$ and $v_i$, respectively. Contracting $e_i$ sends every embedding of $G_i$ to an embedding of $G_{i-1}$, and by induction, the latter is unique. In a disk neighborhood of the embedded edge $e_i$, the attachment points belonging to $A_i$ form one consecutive boundary block, and those belonging to $B_i$ form the other. After contraction, this boundary order becomes the cyclic order around $w_i$. Since $G_{i-1}$ is uniquely embedded, its cyclic order is fixed up to reversal; hence the prescribed splitting must preserve that order, with $A_i$ and $B_i$ consecutive. As a result, the facial closed walk set of embeddings of $G_i$ is identical, by uniquely lifted from the facial closed walk set of $G_{i-1}$, and thus $G_i$ is uniquely embedded too by Lemma \ref{2cell}. Consequently, $T=G_t$ is a uniquely embedded tree expansion of $G_0$.

We still use $\Psi(T)$ to denote the unique embedding of $T$ in $\mathbb{RP}^2$ for every tree expansion $T$. In the end, we prove the if and only condition of $\mathcal{S}^*(k)$ with positive tree expansion, thus completing Theorem \ref{essence2}. By Lemma \ref{key}, we have shown that every edge of $(G,\sigma)\in \mathcal{S}^*(k)$ is deletion-critical in the sense of representativity of its unique embedding, and this condition is also valid to distinguish signed graphs in $\mathcal{S}^*(k)$ from other positive tree expansions of $\mathcal{T}_k$.

Let $(G,\sigma)$ be a positive tree expansion with a canonical embedding such that for every $e\in E(G)$, $\rho(\Psi(G\setminus e))<k$. If there exists an edge cut $\partial_G(X)$ with $|\partial_G^+(X)|<|\partial_G^-(X)|$, then we choose such a cut containing the fewest positive edges. It can be immediately deduced that $\partial_G(X)$ is a minimal cut, and thus both $G[X]$ and $G[X^c]$ are connected. From the cut-curve correspondence stated in Subsection 3.1, $\partial_G(X)$ can be represented by a simple closed curve, and it leaves a single simple arc $\alpha$ in the disk. From the cubicity of $G$ and the premise that $|\partial_G^+(X)|$ is minimum, the positive edges in $\partial_G^+(X)$ form a matching, so after executing the cut-paste procedure with $\alpha$, we obtain another canonical embedding with at most $k-1$ negative edges, then there exists an essential curve lying in the crosscap intersecting $G$ at no more than $k-1$ points, but this is impossible because tree expansion does not decrease representativity. Hence $\sigma$ is a minimum signature of $G$, and $l(G,\sigma)=k$.

On the other hand, arbitrarily select an edge $e\in E(G)$, the condition $\rho(\Psi(G\setminus e))<k$ implies $e$ is traversed by a minimal essential curve $\gamma_e$. From our statement in the proof of Lemma \ref{key}, $\gamma_e$ leaves a simple arc set $\{\alpha_1,\ldots,\alpha_s\}$ in the plane, each of which corresponds to an edge cut $\partial_{G}(X_i)$ by constituting the positive edge set. We assert every $\partial_{G}(X_i)$ is an equilibrated cut, in particular the one containing $e$. On the contrary, if one of them, say $\partial_{G}(X_j)$ satisfies $|\partial_G^+(X_j)|>|\partial_G^-(X_j)|$, then $\gamma_e$ can be transformed into an essential curve $\gamma_e'$ by replacing $\alpha_j$ with a detour through the crosscap, as we have done in the proof of Lemma \ref{key}, with $|\gamma_e'\cap \Psi(G)|\leq k-1$. From the arbitrariness of $e$, the deletion criticality of representativity implies that $(G,\sigma)$ is critically $k$-frustrated, and thus $(G,\sigma)\in \mathcal{S}^*(k)$.
\end{proof}

At last, we give the uniform upper and lower bounds of size of $(G,\sigma)\in \mathcal{S}^*(k)$, hence the finiteness of $\mathcal{S}^*(k)$ holds for every positive integer $k$.

\begin{proof}[\textbf{Proof of Corollary \ref{quad}:}]
Let $(G_0,\sigma)$ be a signed graph obtained from $(G,\sigma)$ by the contractions in Theorem \ref{essence1}, then $|G_0|=|MG_k|=2k^2-k$ because $\Delta Y$ and $Y\Delta$ operations do not change the size of a graph. Every vertex $v$ of degree at least $4$ in $G_0$ is split into cubic vertices, producing $deg_{G_0}(v)-3$ new edges. This gives the upper bound of $|E(G)|$:
\begin{align*}
|E(G)|&=m(G_0)+\sum\limits_{v\in V(G_0)}deg_{G_0}(v)-3=3m(G_0)-3n(G_0)\\
&\leq 3m(G_0)-(m(G_0)+3)=4k^2-2k-3.
\end{align*}
The inequality comes from Euler's formula of $\mathbb{RP}^2$.
\end{proof}

\section{Acknowledgements}
The author was in discussion with Reza Naserasr. The author appreciates the careful guidance from Thomas Zaslavsky on revising the manuscript. Moreover, the author is grateful to Deping Song for pointing out a mistake in a previous version of the proof. The author would also like to thank ChatGPT 5.6 Thinking to inform me of medial graph theory in \cite{medial}.

\section*{Declaration of generative AI and AI-assisted technologies in the manuscript preparation process}
ChatGPT 5.6 Thinking was used in early searching and exploratory discussions on Theorem \ref{essence1}. The AI tool suggested the potential relevance of reference \cite{medial} and provided a lengthy, unintelligible proof which was not adopted. The author then independently studied reference \cite{medial}, and constructed and verified the present proof. The author takes full responsibility for the content of this paper.

\end{document}